\documentclass[12pt,leqno]{amsart}

\usepackage[a4paper,margin=25mm]{geometry}
\usepackage{amsmath,amssymb}

\newtheorem{theorem}{Theorem}
\newtheorem{proposition}[theorem]{Proposition}
\newtheorem{corollary}[theorem]{Corollary}
\theoremstyle{definition}
\newtheorem{example}{Example}

\theoremstyle{remark}
\newtheorem{remark}[theorem]{Remark}

\theoremstyle{definition}

\theoremstyle{plain}
\newtheorem{lemma}[theorem]{Lemma}

\newcommand{\D}{\mathbb D}
\newcommand{\C}{\mathbb C}
\newcommand{\R}{\mathbb R}

\begin{document}
\title{Convexity properties of Kobayashi balls in Lempert domains}
\keywords{Lempert domain, Kobayashi ball, Indicatrix, ($\mathbb C$)-convexity, uniqueness
property of right inverses, complex geodesic}
\thanks{The first named author was partially supported by the Bulgarian National
Science Fund, Ministry of Education and Science of Bulgaria under contract KP-06-N82/6.}
\thanks{The second named author was partially supported by the National Science Centre,
Poland under the Weave UNISONO no. UMO-2025/07/Y/ST1/00146}
\subjclass[2020]{Primary: 32F45, 32F17}
\author{Nikolai Nikolov}
\address{N. Nikolov\\Institute of Mathematics and Informatics\\Bulgarian Academy
of Sciences\\Acad. G. Bonchev 8, 1113 Sofia, Bulgaria \vspace{1mm}
\newline Faculty of Information Sciences\\
State University of Library Studies and Information Technologies\\
Shipchenski prohod 69A, 1574 Sofia, Bulgaria}
\email{nik@math.bas.bg}
\author{W\l odzimierz Zwonek}
\address{Institute of Mathematics, Faculty of Mathematics and Computer Science,
Jagiellonian University, \L ojasiewicza 6, 30--348 Krak\'ow, Poland}
\email{wlodzimierz.zwonek@uj.edu.pl}

\begin{abstract} We clarify relations between different notions of convexity ($\mathbb
C$-convexity, $\mathbb C$-strict convexity) of Lempert domains and their Kobayashi balls.
\end{abstract}

\maketitle

\section{Introduction}

\subsection{General setting of the problem}
We shall be working below in the setting of Lempert domains with the particular stress
put on convex and $\mathbb C$-convex situation. The aim of the paper is to clarify the
relations between different convexity notions of the domain itself and its Kobayashi
balls. In particular, we answer (sometimes partially) problems posed by Vigu\'e and
Znamenski.

First recall that the domain $D\subset\mathbb C^n$ is called {\it Lempert} if the domain
$D$ is taut and the Lempert theorem holds on $D$, i. e. $l_D\equiv c_D$. In this
situation $l_D\equiv k_D\equiv c_D$ and $\kappa_D\equiv\gamma_D$. Therefore, without
loss of generality we may work in this context with the Kobayashi distance $k_D$ and its
infinitesimal version, i. e. Kobayashi-Royden pseudometric $\kappa_D$.

In his fundamental papers Lempert showed his theorem for bounded convex domains and
strongly linearly convex ones. The latter class of domains constitutes a subclass of
another kind of convexity, so called $\mathbb C$-convex domains. Though all bounded
smooth $\mathbb C$-convex domains can be exhausted by strongly linearly convex ones and
thus are Lempert domains it is not known whether all taut $\mathbb C$-convex domains are
Lempert.

For $a\in D$ and $r>0$ put
$$
B_D(a,r):=\{z\in D:k_D(a,z)<r\},\qquad
I_D(a):=\{v\in\mathbb C^n:\kappa_D(a;v)<1\}.$$

It is a direct consequence of the Lempert theorem that the {\it indicatrix} $I_D(a)$ is
a balanced convex domain for all domains satisfying the Lempert theorem (it is
additionally bounded in the case $D$ is Lempert). It is also known that the balls
$B_D(a,r)$ are convex if $D$ is convex (see, e.g., \cite[Lemma~5.4]{Frankel 1989}). It
was posed as an open problem whether the sublevel sets (w.r.t. the Kobayashi distance)
are all $\mathbb C$-convex (see \cite[Problem~3]{Zna 1994}). We solve the problem
positively for domains that may be exhausted by strongly linearly convex ones (in
particular, all bounded smooth $\mathbb C$-convex domains) leaving the open problem that
would be solved if the exhaustion of all taut $\mathbb C$-convex domains would be true.

Our further interest in the paper will be the study of different notions of (strict)
convexity of the Kobayashi balls and (to a smaller extent) indicatrices in Lempert
domains; mainly its relation with the (strict) convexity of the domain itself.

We first recall some notions and properties of Lempert domains.

A holomorphic disc $f:\D\to D$ is called a {\it complex geodesic } if it admits a
holomorphic {\it left inverse} $F:D\to\D$ such that $F\circ f=\operatorname{id}_{\D}$.
The function $F$ is called {\it left inverse of $f$}.

A map $\gamma:(-1,1)\to D$ is called a {\it real geodesic} if
$k_D(\gamma(s),\gamma(t))=p(s,t)$ for $s,t\in(-1,1)$, where $p$ is the Poincar\'e
distance on $\D$.
In Lempert domains every real geodesic segment extends to a real geodesic
(see \cite[Remark~2]{CZ}).

By Lempert's theorem, any two distinct points in Lempert domain lie on a complex
geodesic, and the restriction of any complex geodesic to $(-1,1)$ is a real geodesic.

For every $a\in D$ and $v\in\C^n$ with $\kappa_D(a;v)=1$, there is also a complex
geodesic $f$ with $f(0)=a$ and $f'(0)=v$.

In the original papers of Lempert it is proven that strongly linearly convex domains and
convex domains biholomorphic to bounded domains are Lempert domains. It is also
well-known that Lempert domains are Kobayashi hyperbolic and Kobayashi complete. For these
facts and other properties of invariant functions that we shall use, see
\cite{Lem 1981,L,Lem 1984,Jar-Pfl 2013}.

\subsection{Uniqueness property of right inverses in Lempert domains}
We say that the Lempert domain $D$ has {\it the uniqueness property of right inverses at
$a\in D$} if every holomorphic map $F:D\to\D$ is a left inverse of at most one complex
geodesic $f$ with $f(0)=a$.

As the starting point of the problems studied in the paper may be seen the following
characterization of the uniqueness property of right inverses (see \cite[Proposition~3
and Corollary~5]{CZ}; throughout we refer to the arXiv version).

\begin{proposition}\label{proposition-uniqueness-right-inverses}
Let $D\subset\C^n$ be a Lempert domain and let $a\in D$. Then the following conditions
are equivalent:
\begin{enumerate}
\item[(i)] $D$ has the uniqueness property of right inverses at $a$;
\item[(ii)] for every $b\in D\setminus\{a\}$, there is a unique real
Kobayashi geodesic through $a$ and $b$  (up to reparametrization);
\item[(iii)] every real Kobayashi geodesic through $a$ is induced by a
complex geodesic;
\item[(iv)] $I_D(a)$ is strictly convex and, for every $v\in\C^n$
with $\kappa_D(a;v)=1$, there is a unique complex geodesic $f$ with $f(0)=a$ and
$f'(0)=v$.
\end{enumerate}
\end{proposition}
\begin{proof}
By \cite[Corollary~5]{CZ}, (i), (iii) and (iv) are equivalent. Clearly, (ii) implies
(iii). Conversely, under (i) and (iii), two real geodesics through $a,b$ are induced by
complex geodesics $f,g$ with $f(0)=g(0)=a$ and $f(r)=g(r)=b$ for some $r\in(0,1)$. By
the Schwarz lemma they share a left inverse, so $f=g$.
\end{proof}

\begin{remark}
    Note that by Remark 2 in \cite{CZ} the property (iii) above may be replaced by the uniqueness of geodesic intervals (as defined in \cite{CZ}) with endpoints at $a$ and $b$. 
\end{remark}
In (iv) the uniqueness of complex geodesics follows from strict convexity of the
indicatrix (the idea for the proof below came from a discussion with AI).

\begin{proposition}
Let $D$ be a Lempert domain and let $a\in D$. Then each of conditions (i)--(iv) of
Proposition~\ref{proposition-uniqueness-right-inverses} is equivalent to:
\begin{enumerate}
\item[(v)] $I_D(a)$ is strictly convex.
\end{enumerate}
\end{proposition}
\begin{proof}
Assume that $I_D(a)$ is strictly convex. Let $f,g$ be holomorphic right inverses of
$F\in\mathcal O(D,\mathbb D)$ with $f(0)=g(0)=a$. By \cite[Proposition~4]{CZ}, strict
convexity gives $f'(0)=g'(0)=X\in\partial I_D(a)$.

If $f\ne g$ then we have $f(\zeta)-g(\zeta)=b\zeta^m+O(\zeta^{m+1})$ for some $m\geq2$
and $b\ne0$. For any holomorphic $G:D\to\mathbb D$ with $G(a)=0$ and any $0<r<1$, set
$A(\zeta)=G(f(r\zeta))$ and $B(\zeta)=G(g(r\zeta))$. The functions
\[
H_\pm(\zeta)=\frac{A(\zeta)+B(\zeta)}2
\pm\frac{\zeta^{m-1}+\zeta^{1-m}}4\bigl(A(\zeta)-B(\zeta)\bigr)
\]
are holomorphic in a neighborhood of $\overline{\D}$ except possibly at $0$. Since
$A-B=O(\zeta^m)$, the singularity is removable and $H_\pm(0)=0$. For $|\zeta|=1$ their
values are convex combinations of $A(\zeta)$ and $B(\zeta)$, so the maximum modulus
principle gives $H_\pm(\D)\subset\D$. The Schwarz lemma thus gives
\[
\bigl|G^{\prime}(a)(rX\pm r^m b/4)\bigr|=|H_\pm'(0)|\leq1.
\]
Letting $r\uparrow1$ and taking the supremum over $G$ gives $\kappa_D(a;X\pm
b/4)=\gamma_D(a;X\pm b/4)\leq1$. Consequently, the points $X\pm b/4$ belong to
$\overline{I_D(a)}$. But their midpoint $X$ lies in $\partial I_D(a)$, a contradiction
with the strict convexity of $I_D(a)$.
\end{proof}

\begin{remark}
The argument above also applies to bounded convex domains in reflexive complex Banach
spaces and answers both questions in \cite[pp.~306--307]{Vig 2001}.
\end{remark}

Note that in the uniqueness property of right inverses at a given point we may generally
restrict ourselves to Lempert left inverses.

\begin{lemma}\label{lemma:equivalence-right-inverses} Let $D$ be a Lempert domain with
Lempert left inverses for any geodesics and fix $p\in D$. Then $D$ has the property of
uniqueness of right inverses at $p$ if and only if each Lempert left inverse has at most
one right inverse $f$ with $f(0)=p$.
\end{lemma}
\begin{proof}
To prove the non-trivial implication assume that a mapping $G:D\to \mathbb D$ has two
different complex geodesics $f,g:\mathbb D\to D$ with $f(0)=g(0)=p$,
$G(f(\lambda))=G(g(\lambda))=\lambda$, $\lambda\in\mathbb D$. Then
$\kappa_D(p;f^{\prime}(0))=\kappa_D(p;g^{\prime}(0))=1$. As in our case
$\kappa_D(p;\cdot)$ is a norm we get
\begin{equation}
1\geq\kappa_D\left(p;\frac{f^{\prime}(0)+g^{\prime}(0)}{2}\right)\geq \kappa_{\mathbb
D}\left(G(p);G^{\prime}(p)\left(\frac{f^{\prime}(0)+g^{\prime}(0)}{2}\right)\right)
=\kappa_{\mathbb
D}(0;1)=1,
\end{equation}
so $\kappa_D\left(p;\frac{f^{\prime}(0)+g^{\prime}(0)}{2}\right)=1$. Let $F$ be a
Lempert left inverse for the pair $(p;(f'(0)+g'(0))/2)$. Then $((F\circ f)'(0)+(F\circ
g)'(0))/2=1$, so the Schwarz lemma gives $F\circ f=F\circ g=\operatorname{id}_{\D}$.
\end{proof}

In our paper we relate properties formulated above with the properties of sublevel sets
of important classes of Lempert domains (convex and $\mathbb C$-convex). In order to
formulate the results we introduce some important subclasses of Lempert domains.

\subsection{Different subclasses of Lempert domains} Denote by $\mathcal L$ the class of
Lempert domains. Denote also by $\mathcal L_l$ the class of {\it Lempert domains with
Lempert left inverses}, i. e. Lempert domains such that arbitrary complex geodesic
admits {\it Lempert left inverse}, that is a left inverse $F:D\to\mathbb D$ such that
for any $\lambda\in\mathbb D$ the set $F^{-1}(\lambda)\cap D$ is the intersection of an
affine hyperplane $H_{\lambda}$ with $D$) -- the class of domains studied in a recent
paper \cite{Zwo 2026}.

Recall that bounded smooth $\mathbb C$-convex domains may be exhausted by strongly
linearly convex ones (see \cite{Jac 2006} for the $C^2$-smooth case and
\cite{Edi 2026} for the $C^{1,\alpha}$ case with $\alpha>1/2$) and thus they
belong to $\mathcal L_l$. Denote by $\mathcal L_{slc}$ (respectively, $\mathcal L_{sc}$)
the taut exhaustion of strongly linearly convex (respectively, strongly convex domains).
The following relations hold:
\begin{equation}
\mathcal L_{sc}=\{\text{taut convex domains}\},\;
\mathcal L_{slc}\subset\{\text{taut $\mathbb C$-convex domains}\}.
\end{equation}

It is an open problem whether the last inclusion becomes an equality (see
\cite[Problem~1]{Aiz 1994}). Certainly, $\mathcal L_{slc}\subset \mathcal L_l$. Recall
once more that in the case of bounded convex domains all the Kobayashi balls (and the
indicatrices) are convex. In the case of bounded $\mathbb C$-convex domains (or domains
from the class $\mathcal L_l$) the situation is not so clear. However, we may say
something positive for the class $\mathcal L_{slc}$, which also gives a partial positive
answer to \cite[Problem~3]{Zna 1994}.

In the case of arbitrary Lempert domain we have the following.

\begin{proposition}\label{proposition:lempert-domain}
Let $D$ be a Lempert domain, $p\in D$, $R>0$. Then we have the equality
$l_{B_D(p,R)}(p,z)=c_{B_D(p,R)}(p,z)$ for all $z\in B_D(p,R)$.
\end{proposition}
\begin{proof}
Let $f:\D\to D$ be a complex geodesic through $p$ and $z$ with $f(0)=p$ and left inverse
$F$. For $r\in(0,1)$ with $R=p(0,r)$, the map $f_r:\D\to B_D(p,R)$,
$f_r(\lambda)=f(r\lambda)$, has left inverse $F/r$ restricted to $B_D(p,R)$, by
contractivity.
\end{proof}

\begin{remark}
Can one show that $B_D(p,R)$ under the assumptions of
Proposition~\ref{proposition:lempert-domain} is a Lempert domain?
\end{remark}

\subsection{Different notions of (strict) convexity} As we saw in the preceding sections
different notions of convexity play essential role in the Lempert theory. Except for the
standard convexity a crucial role is played by (geometric) notion of strict convexity:
the convex domain $D\subset\mathbb C^n$ is {\it strictly convex } if any interior point
of a segment with ends in $\overline{D}$ lies in $D$.

We introduce a complex version of the convexity that we shall work with and that is
important in the Lempert theory. A domain $D\subset\mathbb C^n$ is called {\it $\mathbb
C$-convex} if the intersection of $D$ with any complex affine line meeting $D$ is
connected and simply connected. We already saw that a subclass of $\mathbb C$-convex
domains contains Lempert domains (domains from $\mathcal L_{slc}$). We shall consider
below the counterparts of strict convexity in this class of domains. For properties of
$\mathbb C$-convex domains, see, e.g., \cite{And-Pas-Sig 2004}.

Let $D$ be a taut $\mathbb C$-convex domain. We may consider the following three
different notions of strict $\mathbb C$-convexity.

\begin{itemize}
\item for any $z\in\mathbb C^n\setminus D$ any complex affine line passing through $z$
and omitting $D$ may intersect $\overline{D}$ at most one point,
\item for any $z\in \partial D$ there is an affine complex hyperplane $H$ such that
$H\cap\overline{D}=\{z\}$,
\item for any $z\in\partial D$ and any affine supporting complex hyperplane $H$ at $z$ we
have that $H\cap\overline{D}=\{z\}$.
\end{itemize}

Recall that the first and third properties are equivalent (see \cite[Remark~2]{Pfl-Zwo
2019}). The $\mathbb C$-convex domains satisfying this property are called {\it $\mathbb
C$-strictly $\mathbb C$-convex} (or $\mathbb C$-strictly convex if the domain is
convex). The domain satisfying the second property is called {\it $\mathbb C$-extreme
$\mathbb C$-convex}.

If the domain $D$ is $\mathbb C$-strictly $\mathbb C$-convex then it is extreme $\mathbb
C$-convex. In the case the domain is $C^1$-smooth both notions are equivalent. Note also
that strongly linearly convex domains are $\mathbb C$-strictly $\mathbb C$-convex. The
following example shows that extreme $\mathbb C$-convexity is not equivalent to $\mathbb
C$-strict $\mathbb C$-convexity (even in the class of convex domains)

\begin{example}
The bounded convex domain
\[
D=\{(z,w)\in\C^2:|z|^2+|\Im w|<1,\quad |z|^2+|w|^2<2\}
\]
is $\mathbb C$-extreme convex but not \(\C\)-strictly convex.

Indeed, the complex affine line \(L=\{z=1\}\) satisfies
\[
L\cap D=\varnothing,\qquad
L\cap\overline D=\{(1,t):t\in[-1,1]\}.
\]

To show that \(D\) is \(\C\)-extreme convex, fix \(p=(a,b)\in\partial D\).

If \(|a|^2+|\Im b|=1\), choose \(\varepsilon\in\{-1,1\}\) such that \(\varepsilon\Im
b=|\Im b|\), and set
\[
H_p=\{(z,w)\in\C^2:2\bar a(z-a)-i\varepsilon(w-b)=0\}.
\]
For \((z,w)\in H_p\cap\overline D\), taking real parts in the equation of \(H_p\) gives
\[
1\geq |z|^2+\varepsilon\Im w=1+|z-a|^2.
\]
Hence \(z=a\) and then \(w=b\), so \(H_p\cap\overline D=\{p\}\).

Otherwise, \(|a|^2+|b|^2=2\), and the complex hyperplane
\[
H_p=\{(z,w)\in\C^2:\bar a(z-a)+\bar b(w-b)=0\}
\]
meets the closed ball \(\{|z|^2+|w|^2\leq2\}\) only at \(p\). Thus \(H_p\cap\overline
D=\{p\}\) in this case as well.
\end{example}

Recall that any weakly linearly convex balanced domain is convex (see
\cite[Proposition~2]{Nik-Pfl-Zwo 2008}). When assuming the condition of $\mathbb
C$-extreme $\C$-convexity, one gets in this case strict convexity. In particular, in the
class of balanced domains the $\mathbb C$-strict $\mathbb C$-convexity is equivalent to
$\mathbb C$-extremal $\mathbb C$-convexity and usual strict convexity.

\begin{proposition}
Let $D$ be a balanced domain in $\C^n$ that is $\mathbb C$-extreme $\mathbb C$-convex.
Then $D$ is strictly convex.
\end{proposition}
\begin{proof}  By \cite[Proposition~2]{Nik-Pfl-Zwo 2008}, we get that $D$ is convex.
Suppose $D$ is not strictly convex, so there are distinct $x,y\in\partial D$ such that
$a:=\frac{x+y}{2}\in\partial D$. Let $H_a$ be a complex hyperplane with
$\overline{D}\cap H_a=\{a\}$. Let $H_a=\{z\in\C^n:L_a(z)=1\}$, where $L_a$ is a
corresponding linear operator. Then $|L_a(z)|\leq 1$ for all $z\in\overline{D}$.
Actually, if $z\in\overline D$ and $|L_a(z)|>1$, the balanced convexity of $D$ implies
$z/L_a(z)\in D\cap H_a$, a contradiction.

Then
\begin{equation}
1=L_a(a)=\Re L_a(a)=\frac{\Re L_a(x)+\Re L_a(y)}{2}\leq 1.
\end{equation}
As $|L_a(x)|,|L_a(y)|\leq 1$ we get that $L_a(x)=L_a(y)=1$ so $a=x=y$.
\end{proof}

For a discussion on $\C$-strict $\mathbb C$-convexity and related notions, see
\cite{Zim 2017,Pfl-Zwo 2019}.

\section{Kobayashi balls in convex domains}

\subsection{General discussion}
Note that if all Kobayashi balls with a fixed center $a\in D$ are convex, then $D$ is
convex (by exhaustion). Conversely, if $D$ is convex, then all $B_D(a,R)$ and $I_D(a)$
are convex by Lempert's theorem.

We say that a convex domain $D$ is {\it $\C$-proper } if it contains no complex lines.
Recall that any convex domain $D\subset\C^n$ is affinely equivalent to
$G\times\C^{n-k}$, where $0\leq k\leq n$ and $G\subset\C^k$ is a $\C$-proper convex
domain that is additionally Kobayashi hyperbolic and complete (see, e.g., \cite{BS}).

As $B_{G\times\C^{n-k}}((a,b),r)=B_G(a,r)\times\C^{n-k}$ and
$I_D((a,b))=I_G(a)\times\C^{n-k}$, $(a,b)\in G\times\C^{n-k}$ both of them are not
strictly convex if $D\subsetneq\C^n$ is not $\C$-proper.

On the other hand we know that if $D$ is a bounded strictly convex domain, $a\in D$,
then $B_D(a,R)$ and $I_D(a)$ are also strictly convex (see \cite{Vig 2001}). We will
extend this fact.

In the case of tube domain we arrive at an obstruction in producing of $\C$-strictly
convex domains.

\begin{remark}\label{remark:tube-domains}
Let $\Omega\subset\mathbb R^n$ be a convex domain. Define the tube domain with the base
$\Omega$
\begin{equation}
T_{\Omega}:=\{z\in\C^n: \Re(z)\in\Omega\}.
\end{equation}
Certainly tube domains are not strictly convex. For $n\geq 2$ actually more is known;
namely, if $\Omega\neq\R^n$ then it is not $\C$-strictly convex. Actually, if $H$ is a
supporting hyperplane of $\Omega$ at $0\in\partial\Omega$, then $\C H$ is disjoint from
$T_\Omega$ but $\C H\cap\partial T_\Omega=(H\cap\partial\Omega)+iH$.
\end{remark}

\begin{example}
Let $n\geq2$ and let $F_1:\C^{n-1}\to[0,\infty)$ be strictly convex with $F_1(0)=0$. Let
$f:\C\to[0,\infty)$ be convex with $f^{-1}(0)=\{0\}$. Then
\[
D=\bigl\{(z',z_n)\in\C^n:F_1(z')+f(z_n)<1\bigr\}
\]
is a bounded $\C$-strictly convex domain. It is strictly convex if and only if none of
the level sets $\{f(w)=c\}$, $0<c\leq1$, contains a non-trivial segment.
\end{example}
\begin{proof}
Both $F_1$ and $f$ are coercive by convexity and uniqueness of their minima. If
$[p,q]\subset\partial D$, $p\ne q$, then $p'=q'$ and $[p_n,q_n]$ lies in
$\{f=1-F_1(p')\}$. The complex line through $p,q$ contains $(p',0)\in D$. Conversely, a
non-trivial segment in $\{f=c\}$ yields a boundary segment by taking $F_1(z')=1-c$.
\end{proof}

\subsection{Strict convexity of balls in convex domains}
We proceed to the description of $\mathbb C$-strict convexity of Kobayashi balls in
convex domains.

\begin{theorem}\label{pr1}
Let $D\subset\C^n$ be a $\C$-proper convex domain and let $a\in D$. Then the following
conditions are equivalent:
\begin{enumerate}
\item[(i)] $B_D(a,R)$ is $\C$-strictly convex for some $R>0$;
\item[(ii)] $B_D(a,R)$ is strictly convex for every $R>0$;
\item[(iii)] $D$ has the uniqueness property of right inverses at $a$.
\end{enumerate}
\end{theorem}

\begin{proof}
(ii) $\Rightarrow$ (i). This is obvious.

(iii) $\Rightarrow$ (ii). Suppose that $B_D(a,R)$ is not strictly convex for some $R>0$.
Then there are distinct $b,c\in D$ with
$$k_D(a,b)=k_D(a,c)=k_D\!\left(a,\frac{b+c}{2}\right)=R.$$
Let $r\in(0,1)$ be such that $p(0,r)=R$. There are complex geodesics $f,g$ such that
$f(0)=g(0)=a$, $f(r)=b$, and $g(r)=c$. Convexity implies that $h=(f+g)/2$ maps $\D$ into
$D$. Since $k_D(h(0),h(r))=R$, the map $h$ is a complex geodesic. Let $F$ be its left
inverse. By the Schwarz lemma,
$$1=(F\circ h)'(0)=\tfrac12\big((F\circ f)'(0)+(F\circ g)'(0)\big).$$
As $|(F\circ f)^{\prime}(0)|,|(F\circ g)^{\prime}(0)|\leq 1$ the last is not bigger than
$1$ which implies that $F$ is a left inverse of $f$ and $g$, a contradiction.

(i) $\Rightarrow$ (iii). If (iii) fails, we have distinct complex geodesics $f,g$ with
$f(0)=g(0)=a$ and a common left inverse $F$; in particular, $F(f(0))=F(g(0))=0$ and
$F^{\prime}(a)(f^{\prime}(0))=F^{\prime}(a)(g^{\prime}(0))=1$. Applying the Schwarz
lemma to $F\circ f_s$, where $f_s=(1-s)f+sg$, $0\leq s\leq1$, gives
\begin{equation}\label{1}
F\big((1-s)f(\lambda)+sg(\lambda)\big)=\lambda
\qquad(\lambda\in\D,\ 0\leq s\leq1).
\end{equation}
Thus $k_D(a,f_s(\lambda))=p(0,\lambda)$. Moreover, if $f(\lambda)\neq g(\lambda)$ then
by the identity principle $F(z)=\lambda$ for all $z\in l_{\lambda}\cap D$, where
$l_{\lambda}$ is a complex affine line spanned by $f(\lambda)$ and $g(\lambda)$ so then
$k_D(a,z)\geq p(0,\lambda)$ for all $z\in l_{\lambda}\cap D$. For every $R>0$ we may
choose $\lambda$ such that $p(0,\lambda)=R$ and $f(\lambda)\ne g(\lambda)$. But
$[f(\lambda),g(\lambda)]\subset\partial B_D(a,R)$ contradicting the fact that
$l_{\lambda}$ may intersect $\partial B_D(a,R)$ at most one point.
\end{proof}

\begin{proposition}\label{pr2}
Let $D\subsetneq\C^n$ be a  $\C$-strictly convex domain. Then the Kobayashi balls of $D$
are strictly convex.
\end{proposition}

\begin{proof}
First, note that $D$ is $\C$-proper; therefore it is Kobayashi complete and admits a
biholomorphism onto a bounded domain (see, e.g., \cite{BS}). The explicit form of this
biholomorphism implies that every holomorphic disc $f:\D\to D$ has finite radial limits
$f^*(\zeta)$ for a.e. $\zeta\in\partial\D$. Moreover, if $f\ne g$, then $f^*\ne g^*$
a.e.

Suppose that a Kobayashi ball $B_D(p,R)$ is not strictly convex. By Theorem~\ref{pr1},
there are distinct complex geodesics $f,g$ with $f(0)=g(0)=p$ and a common left inverse
$F$. As in earlier proof by convexity and the Schwarz lemma, we obtain the equality
\eqref{1}. For $\lambda\in\D$ with $f(\lambda)\ne g(\lambda)$, put
$L_\lambda=f(\lambda)+\C\big(g(\lambda)-f(\lambda)\big)$. Since $D\cap L_\lambda$ is
connected, the identity theorem and \eqref{1} yield
\begin{equation}\label{2}
F\equiv\lambda\quad\text{on }D\cap L_\lambda.
\end{equation}

As $f\ne g$, we may choose $\zeta\in\partial\D$ such that $A=f^*(\zeta)$ and
$B=g^*(\zeta)$ are finite and distinct. Note that $A,B\in\partial D$. Let $L=A+\C(B-A)$.
If $z=A+c(B-A)\in D$ for some $c\in\C$, then
$$
z_\rho=f(\rho\zeta)+c\big(g(\rho\zeta)-f(\rho\zeta)\big)\longrightarrow z \text{ as
$\rho\to 1$}.
$$
For $\rho$ sufficiently close to $1$, we have $z_\rho\in D$ and $f(\rho\zeta)\ne
g(\rho\zeta)$. Thus \eqref{2} gives $F(z_\rho)=\rho\zeta\to\zeta$, contradicting
$F(z_\rho)\to F(z)\in\D$. Consequently $L\cap D=\varnothing$, a contradiction to $A,B\in
L\cap\partial D$.
\end{proof}

\begin{example}\label{ex2}
Let
$$D=\{z\in\C^2:|z_1|+|z_2|+|z_1-z_2|^2<1\},
\qquad a=(0,0),\quad b=(1/3,1/3).$$
Then $D$ is a bounded convex domain, and for every $r>0$ the ball $B_D(a,r)$ is not
strictly convex, whereas $B_D(b,r)$ is strictly convex.  Moreover, there are no
non-trivial affine discs in the boundary of $D$.
\end{example}

\begin{proof}
Put $\alpha=(\sqrt5-1)/2$. Note that $F(z)=(z_1-z_2)/\alpha$ maps $D$ into $\D$. For
$0\leq t\leq1$, the holomorphic maps $f_t(\lambda)=(t\alpha\lambda,-(1-t)\alpha\lambda)$
take values in $D$ and satisfy $F\circ f_t=\operatorname{id}_{\D}$. Thus, for any $r>0$
and $|\lambda|=\tanh r$, we have $[f_0(\lambda),f_1(\lambda)]\subset \partial B_D(a,r).$

We next show that if $p\neq q$ and $[p,q]\subset\partial D$, then
\begin{equation}\label{3}
|p_1|+|p_2|=\alpha.
\end{equation}
Indeed, equality in the convexity inequality for $|z_1-z_2|^2$ forces
$p_1-p_2=q_1-q_2=s$. Therefore $[p_1,q_1]\subset\{z_1\in\C:|z_1|+|z_1-s|=1-|s|^2\}$.
This set is empty if $1-|s|^2<|s|$, and it is an ellipse if $1-|s|^2>|s|$. Since
$p_1\neq q_1,$ it follows that $1-|s|^2=|s|$, which proves \eqref{3}. The degenerate
ellipse is the segment $[0,s]$. Any affine complex disc in $\partial D$ would thus lie
in a set of the form $\{(w,w-s):w\in[0,s]\}$, which is impossible.

Suppose that some ball centered at $b$ is not strictly convex. By Theorem~\ref{pr1}
there exist distinct complex geodesics $f,g$ with $f(0)=g(0)=b$ and a common left
inverse. As in the proof of Proposition~\ref{pr2}, $f^*,g^*\in\partial D$ and $f^*\ne
g^*$ a.e. By \eqref{1}, $[f^*,g^*]\subset\partial D$. In view of \eqref{3},
$|f_1^*|+|f_2^*|=\alpha$. Subharmonicity yields
$|f_1(0)|+|f_2(0)|\leq\alpha<2/3=|b_1|+|b_2|$, a contradiction.
\end{proof}

\begin{example}\label{ex3}
Let $\Omega\subset\C^n$ be a bounded strictly convex domain and let $H$ be a real open
half-space whose boundary intersects $\Omega$. Then $D=\Omega\cap H$ is not strictly
convex, but all its Kobayashi balls are strictly convex. If $n\geq2$, then $\partial D$
contains non-trivial affine complex discs.
\end{example}
\begin{proof}
After a complex affine change of coordinates, we may assume that $H=\{z\in\C^n:\Re
z_1<0\}$. The relatively open set $\Omega\cap\partial H$ lies in $\partial D$, so $D$ is
not strictly convex and, for $n\geq2$, its boundary contains affine complex discs.
Suppose that some Kobayashi ball centered at $a\in D$ is not strictly convex. By
Theorem~\ref{pr1}, there are distinct complex geodesics $f,g$ with $f(0)=g(0)=a$ and a
common left inverse. By convexity and the Schwarz lemma, $h=(f+g)/2$ is also a complex
geodesic. Their radial limits belong to $\partial D$ almost everywhere. Since $f^*\ne
g^*$ almost everywhere, strict convexity of $\Omega$ gives $h^*\in\Omega$ almost
everywhere. Consequently, $\Re h_1^*=0$ almost everywhere. As $h_1$ is bounded, the
Poisson formula gives $\Re h_1(0)=0$, contradicting $\Re a_1<0$.
\end{proof}

Vigu\'e \cite{Vig 2001} gave an example of a non-strictly convex plane domain with
strictly convex Kobayashi balls.

\begin{example}\label{ex:plane}
A plane domain $D$ is Lempert if and only if it is simply connected and $D\ne\C$. Every
convex domain $D\subsetneq\C$ has strictly convex Kobayashi balls with real-analytic
boundaries, even if $D$ itself is not strictly convex.
\end{example}
\begin{proof}
The first assertion follows from the Riemann mapping theorem and the fact that, by the
identity theorem, any complex geodesic in a plane domain is biholomorphic.

Now let $D\subsetneq\C$ be convex, $a\in D$, and let $f:\D\to D$ be a Riemann map with
$f(0)=a$. For $0<r<1$, we have $B_D(a,p(0,r))=f(r\D)$. Since $f$ is convex,
$$\Re\left(1+\frac{zf''(z)}{f'(z)}\right)>0,\qquad z\in\D.$$
The curvature of the curve $\theta\mapsto f(re^{i\theta})$ is
$$\frac{\Re(1+zf''(z)/f'(z))}{r|f'(z)|}>0,\qquad z=re^{i\theta}.$$
Thus $f(r\D)$ is strictly convex with real-analytic boundary.
\end{proof}

\section{Balls in Lempert ($\mathbb C$-convex) domains}

\subsection{$\mathbb C$-convexity of the balls}
We obtain a similar, though weaker, result for $\mathbb C$-convex domains. Recall that
any $\mathbb C$-convex domain is, up to an affine isomorphism, the Cartesian product
$\mathbb C^k\times D$, where $D$ is a taut $\mathbb C$-convex domain (that is also
Carath\'eodory, Kobayashi complete) (see \cite{Nik-Sar 2007}).

\begin{lemma}\label{lemma:sublevel-sets-strongly-linearly-convex} Let $D$ be a Lempert
linearly convex domain belonging additionally to $\mathcal L_l$ (this is the case for
instance of $C^{1,\alpha}$-smooth bounded $\mathbb C$-convex domains, $\alpha>1/2$) and
let $p\in D$. Then the balls $B_{D}(p,R)$ are linearly convex.

If the domain $D$ is strongly linearly convex then all the Kobayashi balls are smoothly
bounded $\mathbb C$-strictly $\mathbb C$-convex domains.
\end{lemma}
\begin{proof}
Fix $R>0$. Points outside $D$ can be separated from the ball by complex hyperplanes
since $D$ is linearly convex. Choose $z_0$ in $D\setminus B_{D}(p,R)$. Let $f:\D\to D$
be a complex geodesic with $f(0)=p$, $f(t)=z_0$, where $t\in(0,1)$ and $p(0,t)\geq R$.
Let $F:D\to\mathbb D$ be its Lempert left inverse. Then for any $z\in F^{-1}(t)\cap
D=H_t\cap D$ ($H_t$ is an affine hyperplane passing through $z_0$) we get that $F(z)=t$
and hence $k_D(p,z)=c_D(p,z)$ is bigger than or equal to $p(0,t)\geq R$. Consequently,
$H_{t}\cap B_{D}(p,R)=\emptyset$ and finishes the proof in the general case.

Assume now that $D$ is strongly linearly convex. Then $k_D(p,\cdot)$ is smooth off $p$
(see \cite{Lem 1984}). Additionally, the derivative of the function $z\to k_D(p,z)$ is
non-degenerate at any point of $D\setminus\{p\}$ -- this follows from the equality
$k_D(p,f(\lambda))=p(0,|\lambda|)$, $\lambda\in\mathbb D$, where $f:\mathbb D\to D$ is a
complex geodesic with $f(0)=p$. Consequently the ball $B_{D}(p,R)$ is a smooth domain.
Then by results on smooth weakly linearly convex domains (see \cite{And-Pas-Sig 2004})
$B_{D}(p,R)$ is a smooth $\mathbb C$-convex domain.

For $z_0\in\partial B_D(p,R)$ the hyperplane $H_t$ above is the unique complex
supporting hyperplane at $z_0$, by smoothness. If it contains another boundary point
$z_1$, a complex geodesic $g$ with $g(0)=p$, $g(t)=z_1$ has left inverse $F$, by the
Schwarz lemma. This contradicts uniqueness of right inverses (see \cite[Corollary~5 and
Proposition~6]{CZ}; the planar case is immediate).
\end{proof}
\begin{corollary}
Let $D\in\mathcal L_{slc}$. Then its balls are $\C$-convex (and Lempert).
\end{corollary}
\begin{proof}
This follows from Lemma~\ref{lemma:sublevel-sets-strongly-linearly-convex} and stability
of $\C$-convexity and invariant functions under increasing unions. The balls are taut,
being relatively compact $\C$-convex subdomains of the taut $\C$-convex domain $D$ (see
\cite{Nik-Sar 2007}).
\end{proof}
As the bounded $C^{1,\alpha}$-smooth $\mathbb C$-convex domain, $\alpha>1/2$, or the
symmetrized bidisc are exhaustible by strongly linearly convex domains the above applies
to those cases (see \cite{Edi 2026,Pfl-Zwo 2012}).

Recall also that the problem whether Kobayashi balls in $\mathbb C$-convex domains are
$\mathbb C$-convex was formulated in \cite[Problem~3]{Zna 1994}. The above corollary
gives a partial positive answer to that problem.

\begin{remark}
Unlike the case of (strictly) convex domains the result on strongly linearly convex ones
does not seem to be optimal. One could expect the sublevel sets to be strongly linearly
convex in this situation.
\end{remark}

\begin{remark}
The above corollary may be used for instance for the symmetrized bidisc (see
\cite{Pfl-Zwo 2012}) so the balls $B_{\mathbb G_2}(p,R)$ are $\mathbb C$-convex, Lempert
domains. For $p=0$ this gives us effective examples of such domains (the formula for $k_{\mathbb G_2}(0,\cdot)$ is given in
\cite{Agl-You 2001}):
\begin{equation}
\{(s,p)\in\mathbb G_2:2|s-\bar{s}p|+|s^2-4p|<r(4-|s|^2)\},\; r\in(0,1).
\end{equation}
\end{remark}

\subsection{Strict $\mathbb C$-convexity of the balls}
We look for a counterpart of Theorem~\ref{pr1} for domains from the class
$\mathcal L_{slc}$. In the proof below we use
Lemma~\ref{lemma:sublevel-sets-strongly-linearly-convex} and \cite[Remark~2]{Pfl-Zwo
2019}.

\begin{proposition}
Let $D$ be a domain from the class $\mathcal L_{slc}$, $a\in D$.

If $B_D(a,R_0)$ is $\C$-strictly $\C$-convex for some $R_0>0$ then $D$ has uniqueness of
right inverses at $a$.

If $D$ has uniqueness of right inverses at $a$ then $B_D(a,R)$ is $\mathbb C$-extreme
$\mathbb C$-convex for any $R>0$.
\end{proposition}
\begin{proof}
Suppose there are two different complex geodesics $f,g:\mathbb D\to D$ with
$f(0)=g(0)=a$ having the same left inverse $F:D\to \mathbb D$. By
Lemma~\ref{lemma:equivalence-right-inverses} we may assume that $F$ is a Lempert left
inverse. For any $R>0$ we find $\lambda$ with $|\lambda|=r\in(0,1)$ such that
$f(\lambda)\neq g(\lambda)$, where $p(0,r)=R$. But then the points
$f(\lambda),g(\lambda)\in F^{-1}(\lambda)\cap\partial B_D(a,R)$, so the complex
supporting hyperplane containing $F^{-1}(\lambda)$ contains two points lying in the
boundary of $B_D(a,R)$ -- contradiction with the $\mathbb C$-strict $\mathbb
C$-convexity of $B_D(a,R)$.

For the converse, fix $R_0>0$ and let $k_D(a,z_0)=R_0$ and let $f:\mathbb D\to D$ be a
complex geodesic with $f(0)=a$, $f(r_0)=z_0$, $p(0,r_0)=R_0$. Let $F:D\to\mathbb D$ be
Lempert left inverse to $f$. For all $z\in F^{-1}(r_0)\cap D$ we have that $k_D(a,z)\geq
R_0$. If there was some $z\in D$ with $z\neq z_0$, $F(z)=r_0$ and $k_D(a,z)=R_0$, a
complex geodesic $g:\mathbb D\to D$ with $g(0)=a$, $g(r_0)=z$ would satisfy $F\circ
g=\operatorname{id}_{\D}$ by the Schwarz lemma, contradicting uniqueness of right
inverses at $a$. Thus the affine hyperplane $H_{r_0}$ defining $F^{-1}(r_0)$ meets the
closure of $B_D(a,R_0)$ only at $z_0$.
\end{proof}
It is interesting whether in the case of the above proposition we actually have
equivalence of $\mathbb C$-extreme $\mathbb C$-convexity of balls with the $\mathbb
C$-strict $\mathbb C$-convexity.

\end{document}